\documentclass[11pt]{article}
\usepackage{graphicx}
\usepackage{amsmath,amsfonts,amssymb,amsthm}
\title{The typical irregularity of virtual convex bodies}
\author{Rolf Schneider}
\date{}
\newcommand{\Sn}{{\mathbb S}^{n-1}}

\newcommand{\R}{{\mathbb R}}

\newcommand{\cD}{{\mathcal D}}
\newcommand{\K}{{\mathcal K}}

\newcommand{\cP}{{\mathcal P}}
\newcommand{\Rn}{{\mathbb R}^n}

\newcommand{\N}{{\mathbb N}}

\newtheorem{theorem}{Theorem}
\newtheorem{lemma}{Lemma}
\newtheorem{definition}{Definition}

\begin{document}
\maketitle

\begin{abstract}
The semigroup of convex bodies in $\R^n$ with Minkowski addition has a canonical embedding into an abelian group; its elements have been called virtual convex bodies. Geometric interpretations of such virtual convex bodies have been particularly fruitful under the restriction to polytopes. For general convex bodies, mainly the planar case has been studied, as a part of the more general investigation of hedgehogs. Here we restrict ourselves to strictly convex bodies in $\Rn$. A particularly natural geometric interpretation of virtual convex bodies can then be seen in the set of differences of boundary points of two convex bodies with the same outer normal vector. We describe how in the typical case (in the sense of Baire category) this leads to a highly singular object. \\[2mm]
\end{abstract}

\section{Introduction}\label{sec1}

The set $\K^n$ of convex bodies (nonempty, compact, convex subsets) in $\Rn$ ($n\ge 2$), as usual equipped with the vector or Minkowski addition $+$, forms a commutative semigroup with cancellation law. As such, it has a canonical embedding into a commutative group, as follows. On the set of ordered pairs of convex bodies one defines an equivalence relation $\sim$ by $(K,L) \sim (K',L')$ if and only if $K+L'=K'+L$. Let $[K,L]$ denote the equivalence class of $(K,L)$. On the set $G$ of these equivalence classes, an addition $+$ can be defined by $[K,L] + [K',L']= [K+K',L+L']$. Then $(G,+)$ is a commutative group, and the mapping $K\mapsto [K,\{0\}]$ is a semigroup isomorphism of $(\K^n,+)$ into $(G,+)$. The elements of the group $(G,+)$ are called {\em virtual convex bodies}. Virtual polytopes, which are the corresponding objects with $\K^n$ replaced by the set $\cP^n$ of convex polytopes in $\R^n$, were introduced and studied by Pukhlikov and Khovanskii \cite{PK92}. A recent contribution is by Martinez--Maure and Panina \cite{MMP14}. The survey on virtual polytopes by Panina and Streinu \cite{PS15} describes a number of `geometrizations' of virtual polytopes and their applications. A geo\-metrization, or geometric interpretation, of the group of virtual convex bodies is a canonical isomorphism of this group to a group of concrete geometric (or analytic) objects. 

A geometrization of virtual convex bodies that suggests itself is via support functions (for these, see, e.g., Section 1.7 of \cite{Sch14}). We denote by $h_K$ the support function of the convex body $K$, restricted to the unit sphere $\Sn$. With respect to Minkowski addition, support functions have the important property that $h_{K+L}=h_K+h_L$ for $K,L\in\K^n$. Hence, the group of virtual convex bodies is isomorphic to the additive group of differences of support functions of bodies from $\K^n$, via the isomorphism $[K,L]\mapsto h_K-h_L$. The group of differences of support functions is a subgroup of the group of continuous real functions on $\Sn$. The question arises how the differences of support functions can be characterized within the set of continuous functions. For the plane, an answer was given by Martinez--Maure \cite{MM99, MM06}, but the higher-dimensional case seems to be open. 

More intuitively, one might prefer to represent virtual convex bodies not by functions, but by point sets in $\R^n$. A possible starting point is the observation that also support sets of convex bodies behave additively under Minkowski addition. For $u\in \Sn$ we denote by $H_K(u)$ the supporting hyperplane of the convex body $K$ with outer unit normal vector $u$, and by $F_K(u)= H_K(u)\cap K$ the support set of $K$ with normal vector $u$. For $K,L\in \K^n$ we have $H_{K+L}(u)=H_K(u)+H_L(u)$ and $F_{K+L}(u)=F_K(u)+F_L(u)$. Martinez--Maure \cite{MM03} has sketched how the latter fact, together with induction over the dimension, can be used to obtain a geometrization of virtual convex bodies. 

In the following, we restrict ourselves to the subset $\K^n_*\subset \K^n$ of strictly convex bodies. For $K\in\K^n_*$ and $u\in\Sn$, the support set $F_K(u)$ is one-pointed, $F_K(u)= \{x_K(u)\}$ with a unique point $x_K(u)$. It is known that $x_K(u)$ is the gradient of the positively homogeneous support function $h(K,\cdot)$ of $K$ at $u$ (e.g., \cite[Corollary 1.7.3]{Sch14}). Clearly $\{x_K(u): u\in\Sn\} ={\rm bd}\,K$, the boundary of $K$. The mapping $x_K:\Sn \to \R^n$ has been called the {\em reverse spherical image map} (see \cite[p. 88]{Sch14}), or also the {\em reverse Gauss map,} of $K$. It is continuous (see Section \ref{sec2}). 

The mapping
$$ 
[K,L] \mapsto x_K-x_L,\quad K,L\in\K^n_*,
$$
embeds the group $(G_*,+)$ of virtual convex bodies coming from strictly convex bodies into the group of continuous mappings from $\Sn$ into $\Rn$ (with pointwise addition).  We write
$$ x_{K,L} := x_K-x_L$$
and call this the {\em reverse Gauss map} of the pair $(K,L)$. In the following, we study the geometrization of virtual (strictly) convex bodies given by this mapping and its image. 

Generally with a $C^1$ function $h:\Sn\to\R$, Langevin, Levitt and Rosenberg \cite{LLR88} have associated the `hedgehog' $H_h$, defined as the envelope of the family of hyperplanes with equation $\langle x,u\rangle=h(u)$, $u\in\Sn$. They obtain $H_h$ as the image of a unique mapping $x_h:\Sn\to\R^n$ and note that for a point $x_h(u)$ in a smooth part of the (suitably oriented) hypersurface $H_h$, the unit normal vector is $u$. Therefore, the mapping $x_h$ is interpreted as the inverse Gauss map of $H_h$. We remark that hedgehogs have later been thoroughly studied, mainly by Y. Martinez--Maure. In our case, where $h=h_K-h_L$ with $K,L\in\K^n_*$, we shall see in the following (in a more precise form) that smooth parts typically do not exist.

\vspace{-6mm}

\begin{center}
\includegraphics[width=5cm,angle=270]{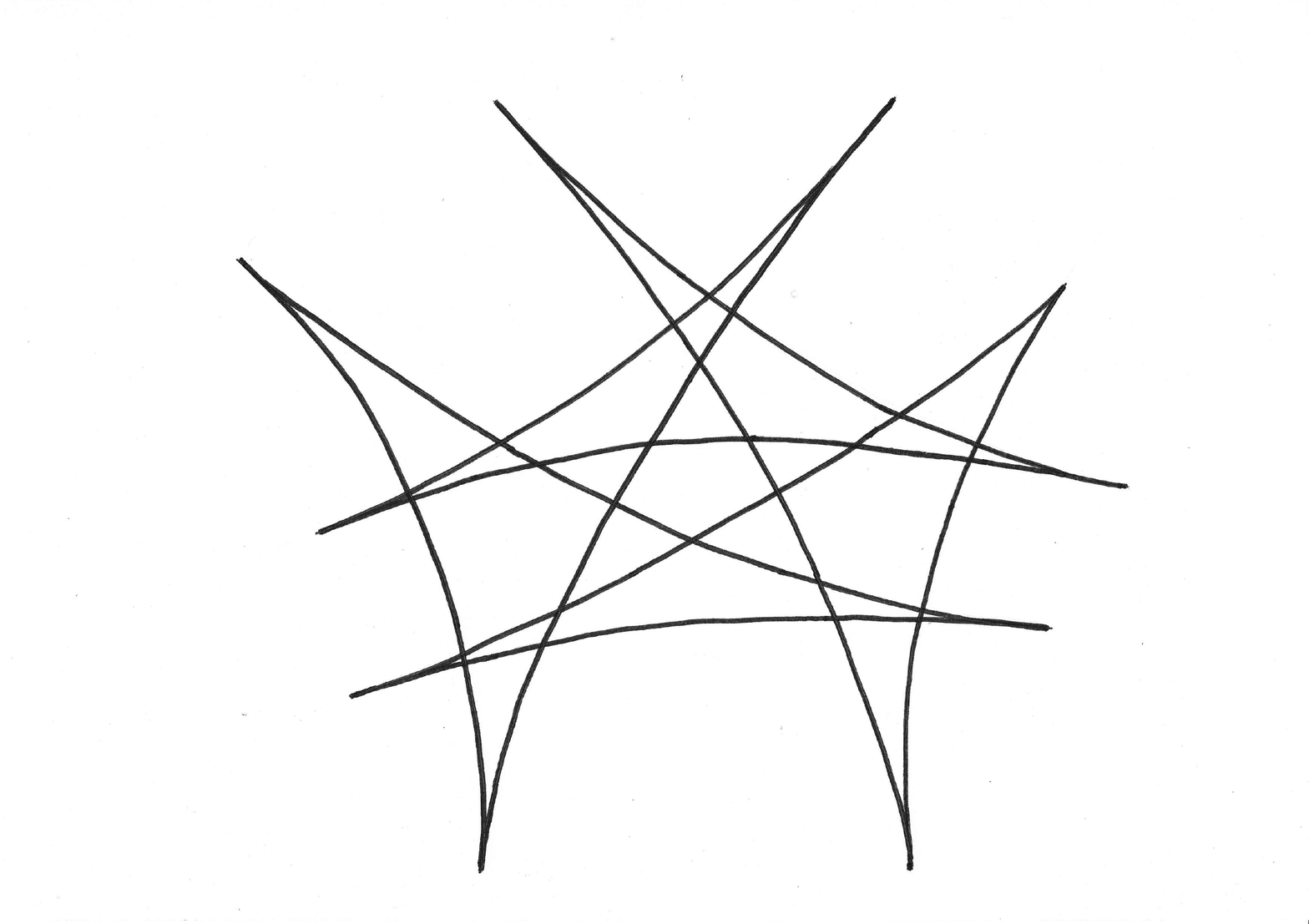}
\end{center}

\vspace*{-5mm}

\noindent {\sc Figure 1.} A sketch, showing how the image set of $x_{K,L}$ may look in simple cases.

\vspace{2mm}

In simple cases, the image set $\cD_{K,L}= x_{K,L}(\Sn)$ looks quite reasonable. For example, let $K,L\in \K^2$ be planar convex bodies of class $C^2$ with positive curvatures and with the property that the radii of curvature of $K$ and $L$, as functions of the outer unit normal vector, have a simple relative behavior. In particular, we assume that they coincide only at finitely many normal vectors $u_1,\dots,u_m$. Then $\cD_{K,L}$ is a closed curve, the points $x_{K,L}(u_1),\dots,x_{K,L}(u_m)$ are cusps, and $\cD_{K,L}$ consists of finitely many convex arcs, each with non-zero radius of curvature except at the endpoints. In particular, for each $u\in{\mathbb S}^1\setminus \{u_1,\dots,u_m\}$, the line $H_K(u)-H_L(u)$ is a tangent line which locally supports $\cD_{K,L}$ at $u$.

The essence of this note is to point out that such nice behavior of the reverse Gauss map $x_{K,L}$ is far from being typical. To make this precise, we understand `typical' in the sense of Baire category. Recall that a topological space $X$ is called a {\em Baire space} if any intersection of countably many dense open subsets of $X$ is dense in $X$. A subset of a Baire space is called {\em comeager} (or {\em residual}\,) if its complement is {\em meager}, which means that it is a countable union of nowhere dense sets. In a Baire space, the intersection of countably many comeager sets is still dense, and for this reason, comeager sets can be considered as `large'. One says that `most' elements of a Baire space $X$ have a property $P$, or that a `typical' element has this property, if the set of all elements with property $P$ is comeager in $X$. The space $\K^n$ of convex bodies, equipped with the Hausdorff metric $\delta$, is a complete metric space and hence a Baire space. It is well known, for example, that most convex bodies in $\K^n$ are smooth and strictly convex (see, e.g., \cite[Theorem 2.7.1]{Sch14}). Surveys on Baire category results in convexity were given by Gruber \cite{Gru85, Gru93} and Zamfirescu \cite{Zam91, Zam09}.) The subspace $\K^n_*$ of strictly convex bodies is a dense $G_\delta$ set in $\K^n$ and hence is also a Baire space. Every set that is comeager in $\K^n_*$ is also comeager in $\K^n$. The product space $(\K^n)^2$, with the metric $d$ defined by $d((K,L),(K',L'))= \delta(K,K')+ \delta(L,L')$ is, of course, also a Baire space, and so is $(K^n_*)^2$. When we show in the following that a set is comeager in $(\K^n_*)^2$, then it is also comeager in $(\K^n)^2$.

The scalar product on $\R^n$ is denoted by $\langle\cdot\,,\cdot\rangle$, the induced norm by $\|\cdot\|$, and $\Sn=\{u\in\R^n:\|u\|=1\}$ is the unit sphere of $\R^n$. Let $u\in \Sn$. In the following, we write
$$ T_u = \{t\in \Sn: \langle t,u\rangle=0\}$$
for the set of unit tangent vectors of $\Sn$ at $u$. For $u\in\Sn$, $t\in T_u$ and $\lambda\in\R$, we use the abbreviation
$$ u \boxplus \lambda t =\frac{u+\lambda t}{\|u+\lambda t\|}.$$

One property of a cusp in the planar example above can be described as follows. Suppose that the mapping $x_{K,L}$ has a cusp at $u\in{\mathbb S}^1$. Then a unit vector $t\in T_u$ can be chosen such that the function $\lambda \mapsto \langle x_{K,L}(u \boxplus \lambda t) - x_{K,L}(u),t\rangle$ is strictly increasing for $\lambda>0$ near $0$ and strictly decreasing for $\lambda <0$ near $0$. In particular, we have 
$$ \langle x_{K,L}(u\boxplus\lambda t)-x_{K,L}(u),t\rangle \cdot\langle x_{K,L}(u\boxplus -\mu t)-x_{K,L}(u),t\rangle >0 $$
for sufficiently small $\lambda, \mu>0$. This is in strong contrast to the behavior of a regular curve with nonzero curvature: if a point on such a curve in the $x,y$-plane moves such that the angle (with the $y$-axis) of the normal vector of the curve at the point moves monotonously in a neighborhood of the direction of the positive $y$-axis, then the $x$-coordinate of the point also moves monotonously. 

The preceding motivates the following definition in higher dimensions.

\begin{definition}\label{D1}
Let $x:\Sn\to\Rn$ be a continuous mapping. Let $u\in\Sn$ and $t\in T_u$. We say that the mapping $x$ 
{\em makes a turn at $u$ in direction $t$} if for every $\varepsilon>0$ there are numbers $0<\lambda,\mu<\varepsilon$ such that
$$ \langle x(u\boxplus \lambda t)-x(u),t\rangle \cdot\langle x(u\boxplus -\mu t)-x(u),t\rangle > 0.$$
Let $\varepsilon>0$. We say that $x$ is {\em $\varepsilon$-tame at $u$ in direction $t$} if 
\begin{equation}\label{1.2}
\langle x(u\boxplus \lambda t)-x(u),t\rangle \cdot\langle x(u\boxplus -\mu t)-x(u),t\rangle \le 0 \quad\mbox{for  } 0<\lambda, \mu<\varepsilon.
\end{equation}
The mapping $x$ is said to be {\em locally tame at $u$ in direction $t$} if $(\ref{1.2})$ holds for some $\varepsilon>0$.
\end{definition}

With this terminology, $x$ makes a turn at $u$ in direction $t$ if and only if it is not locally tame at $u$ in direction $t$.

We note that the reverse Gauss map of a stricly convex body is locally tame everywhere. The same holds for a regular hypersurface of class $C^2$ with nonvanishing curvatures if it has a well-defined inverse Gauss map. 

Now we can formulate our main result. Note that also the unit sphere $\Sn$, and for each $u\in\Sn$ the space $T_u$ of unit tangent vectors at $u$, with their natural topologies, are Baire spaces; the second and third `most' in the following theorem refer to the corresponding space.

\begin{theorem}\label{Thm1}
Most pairs $(K,L)$ of convex bodies in $(\K^n)^2$ have the following properties. The bodies $K$ and $L$ are strictly convex. For most $u\in \Sn$, for most $t\in T_u$, the reverse Gauss map $x_{K,L}$ of $(K,L)$ makes a turn at $u$ in direction $t$.
\end{theorem}

This theorem will be proved in Section \ref{sec3}. The next section contains some preparations.

\section{Auxiliary approximation results}\label{sec2}

On the unit sphere $\Sn$, we use the metric $\Delta$ defined by $\Delta(u,v)=\arccos\langle u,v\rangle$. On the space $\K^n$ of convex bodies, we use the Hausdorff metric, denoted by $\delta$. First, we state a continuity result.

\begin{lemma}\label{L2.1}
The mapping $(K,u)\mapsto x_K(u)$ is continuous on $\K^n_*\times\Sn$.
\end{lemma}

\begin{proof}
Let $((K_i,u_i))_{i\in\N}$ be a sequence in $\K^n_*\times\Sn$ converging to some $(K,u)\in \K^n_*\times\Sn$. The sequence $(x_{K_i}(u_i))_{i\in\N}$ is bounded and hence has a convergent subsequence. Let $x$ be its limit. Since $x_{K_i}(u_i)\in K_i$ for each $i$, we have $x\in K$ (\cite[Thm. 1.8.8]{Sch14}). By the definition of $x_{K_i}(u_i)$, 
$$ \langle x_{K_i}(u_i),u_i \rangle = h(K_i,u_i)\quad \mbox{for }i\in\N,$$
where $h$ denotes the support function, defined on $\K^n\times\Rn$. Since the support function is continuous on the product space (\cite[Lemma 1.8.12]{Sch14}), we get $\langle x,u \rangle = h(K,u)$. Together with $x\in K$ this shows that $x=x_K(u)$, because $K$ is strictly convex. Thus, the sequence $(x_{K_i}(u_i))_{i\in\N}$ has convergent subsequenes, and each such subsequence converges to $x_K(u)$, hence the sequence itself converges to $x_K(u)$. This yields the assertion.
\end{proof}

As a consequence, the mapping $(K,L,u)\mapsto x_{K,L}(u)$ is continuous on $(\K^n_*)^2\times\Sn$.

The proof of Theorem \ref{Thm1} will require some auxiliary results about the approximation of convex bodies by polytopes or by strictly convex bodies. They will be provided in this section.

Let $P\in\K^n$ be a polytope with interior points, and let $F$ be a facet of $P$ (that is, an $(n-1)$-dimensional face). By the {\em normal vector} of $F$ we understand the outer unit normal vector of $P$ at $F$. Suppose that $F_u$ is a facet, with normal vector $u$, and let $F_v$ be a facet, with normal vector $v$, such that $F_u\cap F_v$ is an $(n-2)$-face of $P$. We say in this case that $F_u$ and $F_v$ are {\em adjacent}. The unit vector $t\in T_u$ for which
$$ v= (\cos\alpha)u+(\sin\alpha)t,\quad \alpha=\Delta(u,v),$$
is called the {\em facial tangent vector} of $F_u$ corresponding to $v$. Thus, the facial tangent vectors of $F_u$ are the normal vectors of the facets of $F_u$ relative to the affine hull of $F_u$.

Recall that, for a convex body $K$ and a unit vector $u\in\Sn$, one denotes by $H(K,u)$ the supporting hyperplane of $K$ with outer normal vector $u$ and by $F(K,u)=H(K,u)\cap K$ the face of $K$ with outer normal vector $u$. The closed halfspaces bounded by $H(K,u)$ are denoted by $H^-(K,u)$ and $H^+(K,u)$, where $H^-(K,u)$ is the one with outer normal vector $u$.

For $u\in\Sn$ and $\eta>0$, we write $N(u,\eta)= \{v\in \Sn: \Delta(u,v)<\eta\}$ for the open $\eta$-neighborhood of $u$ in $\Sn$. An {\em $\eta$-net} in $\Sn$ is a finite subset $M\subset\Sn$ with $\bigcup_{u\in M} N(u,\eta)=\Sn$, and an $\eta$-net in $T_u$ is defined similarly.

\begin{lemma}\label{L2.2}
Let $\varepsilon,\eta>0$. Let $M$ be an $\eta$-net in $\Sn$. To every convex body $K\in\K^n$, there exists a polytope $P$ with $\delta(K,P)<\varepsilon$ such that the following holds.\\[1mm]
$\rm (a)$ For $u\in M$, the face $F(P,u)$ is a centrally symmetric facet of $P$.\\[1mm]
$\rm (b)$ For $u\in M$, the facial tangent vectors of $F(P,u)$ form an $\eta$-net in $T_u$.\\[1mm]
$\rm (c)$ Let $u\in M$. If $v\in\Sn$ is such that $F(P,v)$ is a facet of $P$ adjacent to $F(P,u)$, then $\Delta(u,v)<\eta$.\\[1mm]
$\rm (d)$ If $u_1,u_2\in M$, $u_1\not= u_2$, then no facet adjacent to $F(P,u_1)$ is adjacent to $F(P,u_2)$.
\end{lemma}

\begin{proof}
In a first step, we construct a smooth convex body $L$ containing $K$ with $\delta(L,K)<\varepsilon/4$ (for a possible construction, see, e.g., \cite[Sec. 3.4]{Sch14}). That $L$ is smooth means that through each of its boundary points there is a unique supporting hyperplane of $L$. We choose a sequence $(z_i)_{i\in\N}$ that is dense in $\Sn$.  Since
$$ \bigcap_{u\in M} H^-(L,u) \cap \bigcap_{i=1}^k H^-(L,z_i) \to L \quad\mbox{for }k\to\infty,$$
there is a number $k$ such that
$$ Q:= \bigcap_{u\in M} H^-(L,u) \cap \bigcap_{i=1}^k H^-(L,z_i)$$
is a polytope satisfying $\delta(Q,K)<\varepsilon/2$. Since $L$ is smooth, for each $u\in M$ the face $F(Q,u)$ is a facet of $Q$.

For the second step, we choose in $\R^{n-1}$ an $(n-1)$-dimensional, centrally symmetric polytope $C$ with the property that the normal vectors of its $(n-2)$-faces, with respect to $\R^{n-1}$, form an $\eta$-net in $\R^{n-1}\cap\Sn$. Let $u\in M$. We choose a similar image of $C$, denoted by $C_u$, such that $C_u\subset {\rm relint}\,F(Q,u)$, where relint denotes the relative interior. We translate it in direction $u$, obtaining $C_u+\lambda u$, with $\lambda>0$. To each $(n-2)$-face $G$ of $Q$, we choose a unit vector $w_G$ such that the supporting hyperplane $H(Q,w_G)$ satisfies $H(Q,w_G)\cap Q=G$. Then we can choose $\lambda_0>0$ so small that for all $0<\lambda\le\lambda_0$ we have
$$ C_u+\lambda u \subset {\rm int}\,H^-(Q,w_G) \quad\mbox{ for all $(n-2)$-faces $G$ of $Q$}$$
and that the polytope
$$ Q_\lambda:= {\rm conv}\left(Q\cup\bigcup_{u\in M}(C_u+\lambda u)\right)$$
satisfies $\delta(Q_\lambda,K)<\varepsilon$. By the choice of $\lambda_0$, each polytope $C_u+\lambda u$ is a facet of the polytope $Q_\lambda$.

Up to now, the polytope $Q_\lambda$ has, for $0<\lambda\le \lambda_0$, the following properties. It satisfies $\delta(Q_\lambda,K)<\varepsilon$; for each $u\in M$, the face $F(Q_\lambda,u)$ is a centrally symmetric facet of $Q$, and the facial tangent vectors of $F(Q_\lambda,u)$ form an $\varepsilon$-net in $T_u$. 

Let $u\in M$. Let $v\in\Sn$ be such that $F(Q_\lambda,v)$ is a facet of $Q_\lambda$ that is adjacent to $F(Q_\lambda,u)$ and thus intersects it in an $(n-2)$-face $J$. Since $C_u+\lambda u \subset {\rm int}\,H^-(Q,w_G)$ for each $(n-2)$-face of $F(Q,u)$, by construction, the facet $F(Q_\lambda,v)$ is the convex hull of $J$ and some face of $F(Q,u)$. For $\lambda>0$ sufficiently small, the vector $v$ can be made arbitrarily close to $u$. Since $M$ is finite, we can choose $\lambda>0$ so small that for any $u\in M$ and for any normal vector $v$ of a facet of $Q_\lambda$ adjacent to $F(Q_\lambda,u)$, we have $\Delta(u,v)<\eta$. For this $\lambda$, the polytope $P=Q_\lambda$ satisfies conditions (a), (b), (c) of Lemma \ref{L2.2}. It also satisfies (d), since for each $u\in M$, the facets of $Q_\lambda$ adjacent to $F(Q_\lambda,u)$ are contained in the halfspace $H^+(Q,u)$.
\end{proof}

We need an extension of the preceding lemma to simultaneous approximation of two convex bodies. 

\begin{lemma}\label{L2.3}
Let $\varepsilon,\eta>0$. Let $M$ be an $\eta$-net in $\Sn$. To every pair $(K_1,K_2)$ of convex bodies in $\K^n$, there exists a pair $(P_1,P_2)$ of polytopes  with $\delta(K_\nu,P_\nu)<\varepsilon$ for $\nu=1,2$ and having the same normal vectors of facets, such that the following holds.\\[1mm]
$\rm (a)$ For $u\in M$, the face $F(P_\nu,u)$ is a centrally symmetric facet of $P_\nu$, $\nu=1,2$, and $F(P_1,u), F(P_2,u)$ are translates of each other.\\[1mm]
$\rm (b)$ For $u\in M$, the facial tangent vectors of $F(P_\nu,u)$ form an $\eta$-net in $T_u$.\\[1mm]
$\rm (c)$ Let $u\in M$. If $v\in\Sn$ is such that $F(P_1,v)$ is a facet of $P_1$ adjacent to $F(P_1,u)$, then $F(P_2,v)$ is a facet of $P_2$ adjacent to $F(P_2,u)$ (and conversely), and $\Delta(u,v)<\eta$.\\[1mm]
$\rm (d)$ If $u_1,u_2\in M$, $u_1\not= u_2$, then no facet adjacent to $F(P_\nu,u_1)$ is adjacent to $F(P_\nu,u_2)$, $\nu=1,2$.
\end{lemma}

\begin{proof}
We perform the construction in the proof of the preceding lemma for both, $K_1$ and $K_2$. We obtain polytopes $\bar P_1,\bar P_2$ with properties corresponding to those listed in Lemma \ref{L2.2}. Since the similarity factors of the polytopes $C_u$ and the parameter $\lambda>0$ used in the construction can be decreased, if necessary, without altering the success of the construction, we can assume that for each $u\in M$, the facets $F(\bar P_1,u)$ and $F(\bar P_2,u)$ are translates of the same $(n-1)$-polytope $C_u$ and hence of each other. Clearly, the facial tangent vectors of $F(\bar P_\nu,u)$ form an $\eta$-net in $T_u$.

With $0<\mu<1$ we now define
$$ P_1= (1-\mu)\bar P_1 + \mu\bar P_2,\qquad P_2= \mu\bar P_1 + (1-\mu)\bar P_2.$$
We choose $\mu$ so small that $\delta(K_\nu,P_\nu)<\varepsilon$ for $\nu=1,2$. Let $v\in\Sn$ be such that $F(P_1,v)$ is a facet of $P_1$. Then
$$ F(P_1,v)= (1-\mu)F(\bar P_1,v) + \mu F(\bar P_2,v)$$
has dimension $n-1$, hence also $\mu F(\bar P_1,v) + (1-\mu)F(\bar P_2,v)= F(P_2,v)$ has dimension $n-1$, thus $F(P_2,v)$ is a facet of $P_2$ (and conversely). Thus, $P_1$ and $P_2$ have the same normal vectors of facets. 

Let $u\in M$, and suppose that $v\in\Sn$ is such that $F(P_1,v)$ is a facet of $P_1$ adjacent to $F(P_1,u)$. Then $F(P_1,u)\cap F(P_1,v) =:J$ is an $(n-2)$-face of $P_1$ and of $F(P_1,u)$. Since $F(P_2,u)$ is a translate of $F(P_1,u)$, the facet $F(P_2,u)$ has an $(n-2)$-face $J'$ that is a translate of $J$. We have $J'= F(P_2,u)\cap F(P_2,w)$ for some $w\in \Sn$. Since $P_1$ and $P_2$ have the same normal vectors of facets, we must have $w=v$. Here $P_1$ and $P_2$ can be interchanged. 

Moreover, the facets $F(P_\nu,u)$ and $F(\bar P_\nu,u)$, $\nu=1,2$, are all translates of each other, hence each of them has a translate of $J$ as an $(n-2)$-face. Let $v_\nu\in\Sn$ be such that $F(\bar P_v,v_\nu)\cap F(\bar P_\nu,u)$ is a translate of $J$, for $\nu=1,2$. Then the vector $v$ above is equal to either $v_1$ or $v_2$. It follows that $\Delta(u,v)<\eta$.
This completes the proof of the lemma.
\end{proof}

The next lemma concerns approximation by strictly convex bodies. If we approximate a polytope $P$ from the outside by a strictly convex body $K$, then to each facet normal $u$ of $P$, the point $x_K(u)$ is uniquely determined. We need to approximate $P$ in such a way that the position of these points can be prescribed within certain limits.

\begin{lemma}\label{L2.5}
Let $P$ be an $n$-dimensional convex polytope, and let $u_1,\dots,u_k$ be the normal vectors of its facets. For each $i\in\{i,\dots,k\}$, choose a point $z_i\in{\rm relint}\,F(P,u_i)$. Let $\varepsilon_0>0$. There exist a number $0<\varepsilon\le \varepsilon_0$ and a strictly convex body $K$, containing $P$, such that $\delta(K,P)\le\varepsilon_0$ and that 
$$x_K(u_i)= z_i+\varepsilon u_i\quad\mbox{for }i=1,\dots,k.$$
\end{lemma}

\begin{proof}
First we associate with each facet of $P$ a one-parameter family of unbounded closed strictly convex sets, in the following way. Let $F_i=F(P,u_i)$ be a facet of $P$, with normal vector $u_i$. For fixed $i$, we may assume that the given point $z_i$ is equal to $0$ and that ${\rm lin}\,F_i=\R^{n-1}\subset\R^n$. We use polar coordinates in $\R^{n-1}$, representing its points in the form $rv$ with $v\in\Sn\cap \R^{n-1}$ and $r\ge 0$. In $\R^{n-1}$, we choose a stricly convex body $C_i$ containing $F_i$ such that $\delta(F_i,C_i) \le \varepsilon_0/\sqrt{2}$. By $\rho$ we denote the radial function of $C_i$, thus
$$ C_i=\{rv:  v\in \Sn\cap\R^{n-1},\;0\le r\le\rho(v)\}.$$
With a parameter $a>0$, we define
$$ f_a(rv) =a\left[\left(\frac{r}{\rho(v)}\right)^2-1\right],\quad v\in\Sn\cap\R^{n-1},\,r\ge 0.$$
Note that the function $f_a$ attains its minimum at $r=0$, which corresponds to the point $z_i$. Now we define
$$ M(F_i,u_i,a)= \{rv-hu_i: v\in\Sn\cap\R^{n-1},\,r\ge 0,\,h\ge f_a(rv)\}.$$
This set can be viewed as the epigraph of the function $f_a$ (in direction $-u_i$). The nonempty intersections of this set with the translates of $\R^{n-1}$ are homothets of $C_i$, more precisely, for given $h\ge -a$, we have
$$ M(F_i,u_i,a) \cap (\R^{n-1}-hu_i) = \sqrt{(h/a)+1}\,C_i-hu_i.$$
From this and the convexity of the function $f_a$, it follows easily that $M(F_i,u_i,a)$ is convex. To show that $M(F_i,u_i,a)$ is strictly convex, suppose that $H$ is a supporting hyperplane touching the set at two points $r_\nu v_\nu-h_\nu u_i$, $\nu=1,2$. Then $H\cap (\R^{n-1}-h_1u_i)$ is a supporting hyperplane (in dimension $n-1$) of a homothet of $C_i$ at $r_1v_1-h_1u_1$, and $H\cap (\R^{n-1}-h_2u_i)$ is a parallel supporting hyperplane of another homothet of $C_i$ at $r_2v_2-h_2 u_i$. Since $C_i$ is strictly convex, this is only possible if $v_1=v_2$. Then it follows from the strict convexity of the function $f_a$ that $h_1=h_2$. Thus, $M(F_i,u_i,a)$ is strictly convex.

We have associated with each facet $F_i=F(P,u_i)$ the strictly convex set $ M(F_i,u_i,a)$. We can choose $0<\varepsilon\le\varepsilon_0/\sqrt{2}$ so small that 
$$ P\subset M(F_j,u_j,\varepsilon)\quad\mbox{for }j=1,\dots,k.$$
If $z_i+au_i\in M(F_j,u_j,a)$ for some $j\not=i$ and some $a>0$, then also $z_i+a'u_i\in M(F_j,u_j,a')$ for $0<a'<a$. For $a\to 0$, the point $z_i+au_i$ converges to $z_i\in P$. It follows that we can choose $0<\varepsilon\le\varepsilon_0/\sqrt{2}$ so small that in addition
$$ z_i+\varepsilon u_i\in M(F_j,u_j,\varepsilon)\quad\mbox{for }j=1,\dots,k.$$
With this choice of $\varepsilon$, we define
$$ K= \bigcap_{j=1}^k M(F_j,u_j,\varepsilon).$$
Then $K$ is a strictly convex body containing $P$, and for each $i\in\{1,\dots,k\}$ the supporting hyperplane to $K$ with outer normal vector $u_i$ touches $K$ at $z_i+\varepsilon u_i$, that is, $x_K(u_i)=z_i+\varepsilon u_i$.

To estimate the Hausdorff distance of $K$ from $P$, let $x\in K\setminus P$. Then $x$ is separated from $K$ by the affine hull of some facet of $P$, say by $H_i={\rm aff}\,F_i$. Let $x'$ be the image of $x$ under orthogonal projection to $H_i$. If $x'\in F_i$, then $\|x-x'\|\le \varepsilon$. If $x'\notin F_i$, we consider the image $x''$ of $x'$ under orthogonal projection to $C_i$. We have $\|x''-x'\|\le \varepsilon_0/\sqrt{2}$ and hence $\|x-x''\|^2=\|x-x'\|^2+\|x'-x''\|^2\le \varepsilon_0^2$. This shows that $\delta(K,P)\le\varepsilon_0$, which completes the proof.
\end{proof}

\section{Proof of Theorem \ref{Thm1}}\label{sec3}

We prepare the proof of Theorem \ref{Thm1} by some lemmas. For $K,L\in\K^n_*$ and $u\in \Sn$, we define
\begin{equation}\label{3.U} 
U(K,L,u) =  \{ t\in T_u: \mbox{$x_{K,L}$ is locally tame at $u$ in direction $t$}\}
\end{equation}
and, for $k\in\N$,
$$
A_k(K,L,u) =  \{ t\in T_u: \mbox{$x_{K,L}$ is $(1/k)$-tame at $u$ in direction $t$}\}.
$$
Then
\begin{equation}\label{3.3}
U(K,L,u)= \bigcup_{k\in\N} A_k(K,L,u).
\end{equation}

\begin{lemma}\label{L3.1}
Let $(K,L,u)\in(\K^n_*)^2\times\Sn$ and $k\in\N$. The set $A_k(K,L,u)$ is closed in $T_u$.
\end{lemma}

\begin{proof}
Let $(t_i)_{i\in\N}$ be a sequence in $A_k(K,L,u)$ converging to some $t\in T_u$. Let $0<\lambda,\mu < 1/k$. For each $i\in\N$, we have $t_i\in A_k(K,L,u)$ and hence, by Definition \ref{D1},
$$ \langle x_{K,L}(u\boxplus \lambda t_i) - x_{K,L}(u),t_i\rangle \cdot \langle x_{K,L}(u\boxplus -\mu t_i) - x_{K,L} (u),t_i \rangle\le 0.$$
For $i\to\infty$, we have $t_i\to t$ and $u\boxplus\lambda t_i\to u\boxplus\lambda t$; the map $x_{K,L}$ is continuous. Therefore, we obtain
$$ \langle x_{K,L}(u\boxplus\lambda t) - x_{K,L} (u),t\rangle \cdot \langle x_{K,L}(u\boxplus -\mu t) - x_{K,L} (u),t\rangle\le 0.$$
Since this holds for all  $0<\lambda,\mu <1/k$, the map $x_{K,L}$ is $(1/k)$-tame at $u$ in direction $t$, thus $t\in A_K(K,L,u)$. This shows that the latter set is closed.
\end{proof}

For  $K,L\in\K^n_*$, we define
\begin{equation}\label{3.V} 
V(K,L)= \{u\in \Sn: U(K,L,u) \mbox{ is not meager in $T_u$}\}
\end{equation}
and, for $m\in \N$,
\begin{align} \label{3.B}
& B_{k,m}(K,L)\\
&= \{u\in\Sn: \mbox{ There exists $t\in T_u$ such that $N(t,1/m)\cap T_u\subset A_k(K,L,u)$}\}. \nonumber
\end{align}
Then, by (\ref{3.3}),
$$ V(K,L)= \{u\in \Sn: \bigcup_{k\in\N} A_k(K,L,u) \mbox{ is not meager in $T_u$}\}.$$
Suppose that $u\in V(K,L)$. Since $\bigcup_{k\in\N} A_k(K,L,u)$ is not meager in $T_u$, there exists a number $j\in\N$ such that $A_j(K,L,u)$ is not nowhere dense in $T_u$. Since $A_j(K,L,u)$ is closed by Lemma \ref{L3.1}, it has interior points in $T_u$. Therefore, there exists a number $m\in\N$ such that $ u\in B_{j,m}(K,L)$. Conversely, if $ u\in B_{j,m}(K,L)$, then $A_j(K,L,u)$ has interior points in $T_u$, thus $\bigcup_{k\in\N} A_k(K,L,u)$ is not meager, and hence $u\in V(K,L)$. We conclude that
\begin{equation}
\label{3.4} V(K,L)= \bigcup_{k,m\in \N} B_{k,m}(K,L).
\end{equation}

\begin{lemma}\label{L3.2}
Let $K,L\in\K^n_*$ and $k,m\in \N$. The set $B_{k,m}(K,L)$ is closed in $\Sn$.
\end{lemma}

\begin{proof}
Let $(u_i)_{i\in\N}$ be a sequence in $B_{k,m}(K,L)$ converging to some $u\in\Sn$. For each $i\in\N$, we can choose $t_i\in T_{u_i}$ such that 
\begin{equation}\label{3.7}
N(t_i,1/m)\cap T_{u_i} \subset A_k(K,L,u_i).
\end{equation} 
After selecting a subsequence and changing the notation, we can assume that $t_i\to t$ for $i\to\infty$, for some $t\in\Sn$. Then $t\in T_u$. We claim that $N(t,1/m)\cap T_{u} \subset A_k(K,L,u)$.

To show this, let $s\in N(t,1/m)\cap T_u$. Let $0<\lambda,\mu< 1/k$. For each $i\in\N$, we can choose $s_i\in N(t_i,1/m)\cap T_{u_i}$ such that $s_i\to s$ for $i\to\infty$. Then $u_i \boxplus \lambda s_i \to u\boxplus\lambda s$ and $u_i \boxplus -\mu s_i \to u\boxplus -\mu s$ for $i\to\infty$. By (\ref{3.7}), we have $s_i\in A_k(K,L,u_i)$, which implies that 
$$ \langle x_{K,L}(u_i\boxplus \lambda s_i)-x_{K,L}(u_i),s_i\rangle \cdot\langle x_{K,L}(u_i\boxplus -\mu s_i)-x_{K,L}(u_i),s_i\rangle \le 0. $$
For $i\to\infty$, we obtain
$$ \langle x_{K,L}(u\boxplus \lambda s)-x_{K,L}(u),s\rangle \cdot\langle x_{K,L}(u\boxplus -\mu s)-x_{K,L}(u),s\rangle \le 0. $$
Since $0 <\lambda,\mu<1/k$ were arbitrary, this shows that $s\in A_k(K,L,u)$. Since this holds for all $s\in N(t,1/m)\cap T_u$, we have $N(t,1/m)\cap T_u\subset A_k(K,L,u)$. This, in turn, shows that $u\in B_{k,m}(K,L)$. Thus, the latter set is closed.
\end{proof}

Finally, let
\begin{equation}\label{3.M}
{\mathcal M}= \{(K,L)\in (\K^n_*)^2: \mbox{ $V(K,L)$ is not meager in $\Sn$}\}
\end{equation}
and, for $j\in\N$,
\begin{align}\label{3.C} 
&{\mathcal C}_{k,m,j} \\
& =\{(K,L)\in(\K^n_*)^2: \mbox{ There exists $u\in\Sn$ such that $N(u,1/j)\subset B_{k,m}(K,L)$}\}.\nonumber
\end{align}
By (\ref{3.4}) we have
$${\mathcal M}= \{(K,L)\in (\K^n_*)^2: \mbox{ $ \bigcup_{k,m\in \N} B_{k,m}(K,L)$ is not meager in $\Sn$}\}.$$
Suppose that $(K,L)\in{\mathcal M}$. Then there are $k,m\in\N$ such that $B_{k,m}(K,L)$ is not nowhere dense. Since $B_{k,m}(K,L)$ is closed by Lemma \ref{L3.2}, it has interior points. Therefore, there exists a number $j\in\N$ such that $(K,L)\in{\mathcal C}_{k,m,j}$.  We conclude that
\begin{equation}\label{3.10} 
{\mathcal M} \subset \bigcup_{k,m,j\in \N} {\mathcal C}_{k,m,j}.
\end{equation}

\begin{lemma}\label{L3.3}
Let $k,m,j\in \N$. The set ${\mathcal C}_{k,m,j}$ is closed in $(\K^n_*)^2$.
\end{lemma}

\begin{proof}
Let $((K_i,L_i))_{i\in\N}$ be a sequence in ${\mathcal C}_{k,m,j}$ converging to some $(K,L)\in (\K^n_*)^2$. For each $i\in\N$, we can choose $u_i\in\Sn$ such that $N(u_i,1/j)\subset B_{k,m}(K_i,L_i)$. After selecting a suitable subsequence and changing the notation, we can assume that $u_i\to u$ for $i\to\infty$, for some $u\in\Sn$. We claim that $N(u,1/j)\subset B_{k,m}(K,L)$.

To show this, let $z\in N(u,1/j)$. For each $i\in\N$, we can choose $z_i\in N(u_i,1/j)$ such that $z_i\to z$ for $i\to \infty$. Since $z_i\in N(u_i,1/j)\subset B_{k,m}(K_i,L_i)$, we can choose $t_i\in T_{z_i}$ such that 
\begin{equation}\label{3.11}
N(t_i,1/m)\cap T_{z_i}\subset A_k(K_i,L_i,z_i).
\end{equation} After selecting a suitable subsequence and changing the notation, we can assume that $t_i\to t$ for $i\to\infty$, for some $t\in \Sn$. Then $t\in T_z$. We claim that $N(t,1/m)\cap T_z\subset A_k(K,L,z)$. 

To show this, let $s\in N(t,1/m)\cap T_z$. Let $0<\lambda,\mu< 1/k$. For each $i\in\N$, we can choose $s_i\in N(t_i,1/m)\cap T_{z_i}$ such that $s_i\to s$ for $i\to\infty$. Then $z_i \boxplus \lambda s_i \to z\boxplus\lambda s$ and $z_i \boxplus -\mu s_i \to z\boxplus -\mu s$ for $i\to\infty$. By (\ref{3.11}), we have $s_i\in A_k(K_i,L_i,z_i)$, which implies that 
$$ \langle x_{K_i,L_i}(z_i\boxplus \lambda s_i)-x_{K_i,L_i}(z_i),s_i\rangle \cdot\langle x_{K_i,L_i}(z_i\boxplus -\mu s_i)-x_{K_i,L_i}(z_i),s_i\rangle \le 0. $$
Since the mapping $(K,L,u)\mapsto x_{K,L}(u)$ is continuous on $(\K^n_*)^2\times\Sn$ (as follows from Lemma \ref{L2.1}), we conclude with $i\to\infty$ that
$$ \langle x_{K,L}(z\boxplus \lambda s)-x_{K,L}(z),s\rangle \cdot\langle x_{K,L}(z\boxplus -\mu s)-x_{K,L}(z),s\rangle \le 0. $$
Since $0 <\lambda,\mu<1/k$ were arbitrary, this shows that $s\in A_k(K,L,z)$. Since this holds for all $s\in N(t,1/m)\cap T_z$, we have $N(t,1/m)\cap T_z\subset A_k(K,L,z)$. This shows that $z\in B_{k,m}(K,L)$. Since this holds for all  $z\in N(u,1/j)$, we have proved that $N(u,1/j) \subset B_{k,m}(K,L)$. By definition, this means that $(K,L)\in{\mathcal C}_{k,m,j}$. Thus, the latter set is closed.
\end{proof}

\begin{lemma}\label{L3.4}
Let $k,m,j\in \N$. The set ${\mathcal C}_{k,m,j}$ is nowhere dense in $(\K^n_*)^2$.
\end{lemma}

\begin{proof}
Since ${\mathcal C}_{k,m,j}$ is closed in $(\K^n_*)^2$ by Lemma \ref{L3.3}, we have to show that it has empty interior in $(\K^n_*)^2$. For this, let $(\bar K_1,\bar K_2)\in(\K^n_*)^2$ and $\varepsilon_0>0$ be given. We shall show that the $\varepsilon_0$-neighborhood of $(\bar K_1,\bar K_2)$ in $(\K^n_*)^2$ contains an element of the complement of ${\mathcal C}_{k,m,j}$

Choose $0<\eta<\min\{1/m,1/j\}$ and so small that
\begin{equation}\label{3.12}
\Delta(u,u\boxplus \lambda t)<\eta\quad\mbox{implies}\quad \lambda<1/k.
\end{equation}
(Here $u\in\Sn$, $t\in T_u$, and we have $\lambda=\tan\Delta(u,u\boxplus\lambda t)$.) Choose an $\eta$-net $M$ in $\Sn$.

By Lemma \ref{L2.3}, there exists a pair $(P_1,P_2)$ of polytopes with $\delta(\bar K_\nu,P_\nu)<\varepsilon_0/4$ for $\nu=1,2$ and having the same normal vectors of facets, such that the following holds.\\[1mm]
$\rm (a)$ For $u\in M$, the face $F(P_\nu,u)$ is a centrally symmetric facet of $P_\nu$, $\nu=1,2$, and $F(P_1,u), F(P_2,u)$ are translates of each other.\\[1mm]
$\rm (b)$ ´For $u\in M$, the facial tangent vectors of $F(P_\nu,u)$ form an $\eta$-net in $T_u$.\\[1mm]
$\rm (c)$ Let $u\in M$. If $v\in\Sn$ is such that $F(P_1,v)$ is a facet of $P_1$ adjacent to $F(P_1,u)$, then $F(P_2,v)$ is a facet of $P_2$ adjacent to $F(P_2,u)$ (and conversely), and $\Delta(u,v)<\eta$.\\[1mm]
$\rm (d)$ If $u_1,u_2\in M$, $u_1\not= u_2$, then no facet adjacent to $F(P_\nu,u_1)$ is adjacent to $F(P_\nu,u_2)$, $\nu=1,2$.

Let $u\in M$. By property (a) above, there exists a translation vector $y$ such that
$$ F(P_2,u)=F(P_1,u)+y.$$
Let $(t,-t)$ be a pair of facial tangent vectors of $F(P_1,u)$ (which is centrally symmetric) and hence of $F(P_2,u)$. Let $G^+$ be the $(n-2)$-face of $P_1$ that corresponds to the facial tangent vector $t$, and let $G^-$ be the $(n-2)$-face of $P_1$ corresponding to $-t$. Let $v^+\in\Sn$ be the vector such that the facet $F(P_1,v^+)$ intersects $F(P_1,u)$ in $G^+$. Similarly, let $v^-\in \Sn$ be such that $F(P_1,v^-)\cap F(P_1,u)=G^-$. Since $P_1$ and $P_2$ have the same normal vectors of facets, we then have $F(P_2,v^+)\cap F(P_2,u)= G^++y$ and $F(P_2,v^-)\cap F(P_2,u)= G^-+y$. The facets $F(P_1,v^+)$ and $F(P_2,v^+)$ need not be translates of each other, but $G^+$ is an $(n-2)$-face of $F(P_1,v^+)$ and $G^++y$ is an $(n-2)$-face of $F(P_2,v^+)$. Therefore, we can choose points 
$$ p^+,q^+\in{\rm relint}\, F(P_1,v^+) \quad\mbox{such that}\quad q^+= p^+-\beta^+ u+\gamma^+ t\quad \mbox{with }\beta^+,\gamma^+>0$$
and that
$$ p^++y,\,q^++y\in{\rm relint}\,F(P_2,v^+).$$
Similarly, we can choose points 
$$ p^-,q^-\in{\rm relint}\, F(P_1,v^-) \quad\mbox{such that}\quad q^-= p^--\beta^- u-\gamma^- t\quad \mbox{with }\beta^-,\gamma^->0$$
and that
$$ p^-+y,\,q^-+y\in{\rm relint}\,F(P_2,v^-).$$
Now we define
$$ z_{v^+}^1= p^+, \quad z_{v^-}^1= q^-, \quad z_{v^+}^2 = q^++y, \quad z_{v^-}^2= p^-+y.$$
Then we have
\begin{equation}\label{3.13}
z_{v^+}^2= z_{v^+}^1-\beta^+u+\gamma^+t+y,\qquad z_{v^-}^2= z_{v^-}^1+\beta^-u+\gamma^-t+y.
\end{equation}

Since the vector $u\in M$ and the facial tangent vector $t$ of $F(P_\nu,u)$ were arbitrary, we have in this way associated, with each normal vector $v$ of a facet adjacent to $F(P_\nu,u)$ with $u\in M$, a point $z_v^\nu\in{\rm relint}\,F(P_\nu,v)$. Property (d) above assures that these points are well-defined, since $F(P_\nu,v)$ cannot be adjacent to $F(P_\nu,u_1)$ and $F(P_\nu,u_2)$ for different $u_1,u_2\in M$. For the remaining normal vectors $w$ of facets of $P_\nu$, we choose points $z_w^\nu\in{\rm relint}\,F(P_\nu,w)$ arbitrarily, but in such a way that
\begin{equation}\label{3.14} 
z_u^2 = z_u^1+y\quad\mbox{for } u\in M.
\end{equation}

To each of the polytopes $P_1,P_2$ and to the points chosen in the relative interiors of their facets, we now apply Lemma \ref{L2.5}. This yields a number $0< \varepsilon \le \varepsilon_0/4$ and strictly convex bodies $K_\nu$ containing $P_\nu$ such that $\delta(K_\nu,P_\nu) \le\varepsilon_0/4$ and hence $\delta(K_\nu,\bar K_\nu)<\varepsilon/2$, therefore
$d((K_1,K_2),(\bar K_1,\bar K_2)) \le\varepsilon_0$, having, in particular, the following property. If $u\in M$, if $(t,-t)$ is a pair of facial tangent vectors of $F(P_1,u)$, and if $v^+$ and $v^-$ are the normal vectors of the facets of $P_1$ that are adjacent to $F(P_1,u)$ and correspond to $t$ and $-t$, respectively, then 
\begin{equation}\label{3.15}
x_{K_\nu}(v^+)= z_{v^+}^\nu+\varepsilon v^+,\qquad x_{K_\nu}(v^-)= z_{v^-}^\nu+\varepsilon v^-, \qquad x_{K_\nu}(u)= z_u^\nu+\varepsilon u.
\end{equation}
Because of (\ref{3.12}) and $\Delta(u,v^+)<\eta$, $\Delta(u,v^-)<\eta$ (by property (c) above), we can write $v^+= u\boxplus \lambda t$ and $v^-= u\boxplus -\mu t$ with $0<\lambda,\mu <1/k$. From (\ref{3.13}), (\ref{3.14}), (\ref{3.15}) we now get, with positive numbers $\beta^+,\gamma^+$ (depending on $v^+$),
\begin{align*}
& \langle x_{K_1,K_2}(u\boxplus \lambda t)- x_{K_1,K_2}(u),t\rangle\\
&=\langle x_{K_1}(v^+) -x_{K_2}(v^+) -[x_{K_1}(u)-x_{K_2}(u)],t\rangle\\
&= \langle z_{v^+}^1+\varepsilon v^+-(z_{v^+}^2+\varepsilon v^+)+y,t\rangle\\
& =\langle \beta^+u-\gamma^+t,t\rangle = -\gamma^+  <0.
\end{align*}
Similarly, with positive numbers $\beta^-,\gamma^-$
\begin{align*}
& \langle x_{K_1,K_2}(u\boxplus -\mu t)- x_{K_1,K_2}(u),t\rangle\\
&=\langle x_{K_1}(v^-) -x_{K_2}(v^-) -[x_{K_1}(u)-x_{K_2}(u)],t\rangle\\
&= \langle z_{v^-}^1+\varepsilon v^--(z_{v^-}^2+\varepsilon v^-)+y,t\rangle\\
& =\langle -\beta^-u-\gamma^-t,t\rangle = -\gamma^-  <0.
\end{align*}
By Definition \ref{D1}, this means that the reverse Gauss mapping $x_{K_1,K_2}$ is not $(1/k)$-tame at u in direction $t$.  Since the vectors $t$ for which this holds form an $\eta$-net in $T_u$ (by property (b) above), and $\eta<1/m$, we have proved that $u\notin B_{k,m}(K_1,K_2)$, according to definition (\ref{3.B}). This holds for each $u\in M$. Since $M$ is an $\eta$-net in $\Sn$, and $1/j<\eta$, we have shown that $(K_1,K_2)\notin {\mathcal C}_{k,m,j}$, according to definition (\ref{3.C}). Since $d((K_1,K_2),(\bar K_1,\bar K_2)) <\varepsilon_0$, this completes the proof.
\end{proof}

We are now in a position to finish the proof of Theorem \ref{Thm1}. From (\ref{3.10}) and Lemma \ref{L3.4} it follows that ${\mathcal M}$ (defined by (\ref{3.M})) is meager in $(\K^n_*)^2$. Hence, its  complement ${\mathcal M}^c= (\K^n_*)^2\setminus {\mathcal M}$ is comeager. This means that for most pairs $(K,L)\in (\K^n_*)^2)$, the set $V(K,L)$ (defined by (\ref{3.V})) is meager in $\Sn$. That $V(K,L)$ is meager in $\Sn$, implies that the complement $V(K,L)^c=\Sn\setminus V(K,L)$ is comeager in $\Sn$, that is, for most $u\in\Sn$, the set $U(K,L,u)$ (defined by (\ref{3.U})) is meager in $T_u$. That $U(K,L,u)$ is meager in $T_u$, implies that the complement $U(K,L,u)^c= T_u\setminus U(K,L,u)$ is comeager in $T_u$, thus  for most $t\in T_u$, the map $x_{K,L}$ is not locally tame at $u$ in direction $t$. Altogether we have shown that for most pairs $(K,L)\in (\K^n_*)^2$, for most $u\in\Sn$, for most $t\in T_u$, the map $x_{K,L}$ is not locally tame at $u$ in direction $t$, equivalently that it makes a turn at $u$ in direction $t$ . This was the assertion of Theorem \ref{Thm1}.

\vspace{3mm}

\noindent Author's address:\\[2mm]
Rolf Schneider\\
Mathematisches Institut, Albert-Ludwigs-Universit{\"a}t\\
D-79104 Freiburg i. Br., Germany\\
E-mail: rolf.schneider@math.uni-freiburg.de

\end{document}